\documentclass[12pt]{article}
\usepackage{amsfonts,amsmath}
\def \beq {\begin{eqnarray}}
\def \eeq {\end{eqnarray}}
\def \beqn {\begin{eqnarray*}}
\def \eeqn {\end{eqnarray*}}

\newcommand{\halmos}{\rule{1ex}{1.4ex}}

\newcounter{for}[section]

\newtheorem{itlemma}{Lemma}[section]
\newtheorem{itproposition}[itlemma]{Proposition}
\newtheorem{theorem}[itlemma]{Theorem}
\newtheorem{itcorollary}[itlemma]{Corollary}
\newtheorem{itremark}[itlemma]{Remark}
\newtheorem{itremarks}[itlemma]{Remarks}
\newtheorem{itdefinition}[itlemma]{Definition}
\newtheorem{itexample}[itlemma]{Example}

\newenvironment{fact}{\begin{itfact}\rm}{\end{itfact}}
\newenvironment{claim}{\begin{itclaim}\rm}{\end{itclaim}}
\newenvironment{lemma}{\begin{itlemma}}{\end{itlemma}}
\newenvironment{remark}{\begin{itremark}\rm}{\end{itremark}}
\newenvironment{remarks}{\begin{itremarks} \rm}{\end{itremarks}}
\newenvironment{corollary}{\begin{itcorollary}}{\end{itcorollary}}
\newenvironment{proposition}{\begin{itproposition}}{\end{itproposition}}
\newenvironment{definition}{\begin{itdefinition}\rm}{\end{itdefinition}}
\newenvironment{example}{\begin{itexample}\rm}{\end{itexample}}
\newenvironment{proof}{\noindent {\em Proof}.\ \
}{\hspace*{\fill}$\halmos$\medskip}
\newcommand{\be}[1]{\addtocounter{for}{1} \begin{equation}\label{#1}}
\newcommand{\ee}{\end{equation}}
\newcommand{\bl}[1]{\begin{lemma}\label{#1}}
\newcommand{\br}[1]{\begin{remark}\label{#1}}
\newcommand{\brs}[1]{\begin{remarks}\label{#1}}
\newcommand{\bt}[1]{\begin{theorem}\label{#1}}
\newcommand{\bd}[1]{\begin{definition}\label{#1}}
\newcommand{\bp}[1]{\begin{proposition}\label{#1}}
\newcommand{\bc}[1]{\begin{corollary}\label{#1}}
\newcommand{\bfact}[1]{\begin{fact}\label{#1}}
\newcommand{\bex}[1]{\begin{example}\label{#1}}
\newcommand{\ec}{\end{corollary}}
\newcommand{\efact}{\end{fact}}
\newcommand{\eex}{\end{example}}
\newcommand{\el}{\end{lemma}}
\newcommand{\er}{\end{remark}}
\newcommand{\ers}{\end{remarks}}
\newcommand{\et}{\end{theorem}}
\newcommand{\ed}{\end{definition}}
\newcommand{\ep}{\end{proposition}}
\newcommand{\epr}{\end{proof}}
\newcommand{\bpr}{\begin{proof}}
\newcommand{\bcl}[1]{\begin{claim}\label{#1}}
\newcommand{\ecl}{\end{claim}}
\newcommand{\ecs}{\end{corollary}}
\newcommand{\eers}{\end{exercise}}
\newcommand{\eexs}{\end{example}}
\newcommand{\eems}{\end{example}}
\newcommand{\els}{\end{lemma}}
\newcommand{\eles}{\end{lemmaex}}
\newcommand{\ets}{\end{theorem}}
\newcommand{\eds}{\end{definition}}
\newcommand{\eps}{\end{proposition}}

\newcommand{\bi}{\begin{itemize}}
\newcommand{\ei}{\end{itemize}}
\newcommand{\ben}{\begin{enumerate}}
\newcommand{\een}{\end{enumerate}}

\def\vbar{\mathchoice{\vrule height6.3ptdepth-.5ptwidth.8pt\kern-.8pt}
   {\vrule height6.3ptdepth-.5ptwidth.8pt\kern-.8pt}
   {\vrule height4.1ptdepth-.35ptwidth.6pt\kern-.6pt}
   {\vrule height3.1ptdepth-.25ptwidth.5pt\kern-.5pt}}
\def\fudge{\mathchoice{}{}{\mkern.5mu}{\mkern.8mu}}
\def\bbc#1#2{{\rm \mkern#2mu\vbar\mkern-#2mu#1}}
\def\bbb#1{{\rm I\mkern-3.5mu #1}}
\def\bba#1#2{{\rm #1\mkern-#2mu\fudge #1}}
\def\bb#1{{\count4=`#1 \advance\count4by-64 \ifcase\count4\or\bba A{11.5}\or
   \bbb B\or\bbc C{5}\or\bbb D\or\bbb E\or\bbb F \or\bbc G{5}\or\bbb H\or
   \bbb I\or\bbc J{3}\or\bbb K\or\bbb L \or\bbb M\or\bbb N\or\bbc O{5} \or
   \bbb P\or\bbc Q{5}\or\bbb R\or\bbc S{4.2}\or\bba T{10.5}\or\bbc U{5}\or
   \bba V{12}\or\bba W{16.5}\or\bba X{11}\or\bba Y{11.7}\or\bba Z{7.5}\fi}}

\def \T {{\cal{T}}}
\def \L {{\cal{L}}}

\def \H {{\cal{H}}}

\def \D {{\cal{D}}}

\def \S {{\cal{S}}}

\newcommand{\ba}[1]{\addtocounter{for}{1} \begin{eqnarray}\label{#1}}
\newcommand{\ea}{\end{eqnarray}}

\def\sqr#1#2{{\vcenter{\vbox{\hrule height .#2pt
                             \hbox{\vrule width .#2pt height#1pt \kern#1pt
                                   \vrule width .#2pt}
                             \hrule height .#2pt}}}}

\def\pmb#1{\setbox0=\hbox{#1}%
   \kern-.025em\copy0\kern-\wd0
   \kern.05em\copy0\kern-\wd0
 \kern-.025em\raise.0433em\box0 }
\def\sqr#1#2{{\vcenter{\vbox{\hrule height.#2pt
     \hbox{\vrule width.#2pt height#1pt \kern#1pt
   \vrule width.#2pt}\hrule height.#2pt}}}}

\def\B{{\mathcal B}}
\def\E{{\mathcal E}}
\def \RR{{\mathcal R}}

\def \I{{\mathcal I}}
\def\N{{\mathbb N}}
\def\Z{{\mathbb Z}}
\def\R{{\mathbb R}}

\def\P{{\mathbb P}}
\def\TT{{\mathbb T}}
\def\EE{{\mathbb E}}

\def\capa{\text{cap}}
\def\leaf{\text{leaf}}

\def\v{\varphi}

\def\bs{\backslash}
\def\reff#1{(\ref{#1})}

\def \ind {\hbox{ 1\hskip -3pt I}}
\newcommand {\cro}[1] {\left[ {#1} \right]}
\newcommand {\acc}[1] {\left\{ {#1} \right\}}
\newcommand {\pare}[1] {\left( {#1} \right)}

\newcommand {\bra}[1] {\left< {#1} \right>}
\newcommand {\under}[1] {\underline{#1}}

\begin{document}

\title{On large intersection and self-intersection local times in dimension five or more.}
\author{Amine Asselah \\ Universit\'e Paris-Est\\
Laboratoire d'Analyse et de Math\'ematiques Appliqu\'ees\\
UMR CNRS 8050\\ amine.asselah@univ-paris12.fr}
\date{}
\maketitle
\begin{abstract}
We show a remarkable similarity between strategies to realize a 
large intersection or self-intersection local times in dimension five or more.
This leads to the same rate functional for large deviation principles for the
two objects obtained respectively by Chen and M\"orters 
in~\cite{chen-morters}, and by the present author in~\cite{A07}. 
Also, we present a new estimate for the distribution of {\it high}
level sets for a random walk, with application to the geometry of the
intersection set of two {\it high} level sets
of the local times of two independent random walks.
\end{abstract}

{\em Keywords and phrases}: self-intersection local times, large deviations, random walk.

{\em AMS 2000 subject classification numbers}: 60K35, 82C22,
60J25.

{\em Running head}: Intersection local times.

\section{Introduction} \label{sec-intro}
Our purpose, in this paper, is twofold: (i) to unravel some connections
between the rate functionals for the large deviations for 
intersection and self-intersection local times in dimension five or more, and
(ii) to explore the geometry of intersection set of two level sets
of the local times of two independent random walks.

To describe the two quantities we are comparing, a few notations are needed.
Thus, we consider two independent
aperiodic simple random walks on the cubic lattice $\Z^d$, with
$d\ge 5$. More precisely, if $S_n$ is the position of the first walk 
at time $n\in \N$, then
$S_{n+1}$ chooses uniformly at random a site of $\acc{z\in \Z^d: |z-S_n|\le 1}$,
where for $z=(z_1,\dots,z_d)\in \Z^d$, the $l^1$-norm is $|z|:=|z_1|+\dots+|z_d|$.
When $S_0=x$, we denote the law of this walk by $P_x$, and its expectation by $E_x$.
All quantities related to the second walk differ with a tilda, whereas $\P=P_0
\otimes\tilde P_0$. For $n\in \N\cup \acc{\infty}$, we
denote the local times in a time period $[0,n]$ by $\acc{l_n(z),z\in \Z^d}$ with
$l_n(z)=\sum_{k=0}^n\ind\{S_k=z\}$. The intersection local times of two
random walks, in an infinite time horizon, is
\[
\bra{l_\infty,\tilde l_\infty}=\sum_{z\in \Z^d} l_\infty(z)\tilde l_\infty(z),\quad
\text{and}\quad \EE\cro{\bra{l_\infty,\tilde l_\infty}}=\sum_{z\in \Z^d}
G_d(z)^2<\infty,
\]
where the Green's function $G_d$ is square summable in dimension 5 or more.
On the other side, the self-intersection local times, in a time period $[0,n]$,
is
\[
\bra{l_n,l_n}=\sum_{z\in \Z^d} l_n^2(z),\quad
\text{and}\quad \lim_{n\to\infty}\frac{E_0\cro{\bra{l_n,l_n}}}{n}=2G_d(0)-1.
\]

On one hand, Chen and M\"orters in~\cite{chen-morters}
have obtained a large deviations principle for $\{\langle l_\infty,\tilde l_\infty
\rangle\ge t\}$ for $t$ large, in dimension 5 or
more, by an elegant asymptotic estimation of the moments, improving on the
pioneering work of Khanin, Mazel, Shlosman and Sinai in~\cite{KMSS}. Their
method provides a variational formula for the rate functional.

On the other hand, we have obtained a large deviation principle for
$\{\langle l_n,l_n\rangle -\EE[\langle l_n,l_n\rangle]\ge \xi n\}$ for $\xi>0$
by a direct study of the contribution of each level set of the local times.
Our approach provides information on the optimal strategy, but no information on the rate
functional since it eventually relies on a subadditive argument.

Our first result establishes that in spite of a different time-horizon, 
$\langle l_\infty,\tilde l_\infty \rangle$ and $\langle l_n,l_n\rangle$ 
yield the same rate functions (up to a natural factor of 2).

\bp{prop-limit}
Consider a random walk in dimension 5 or more. Then, for any $\xi>0$
\be{LDP-limit}
\lim_{n\to\infty} \frac{1}{\sqrt n} \log P\pare{ \bra{l_n,l_n}-E_0[\bra{l_n,l_n}]>\xi n}
=-\I_{CM}{\sqrt \xi},
\ee
where $\I_{CM}$ is the rate functional for $\langle l_\infty,\tilde l_\infty \rangle$,
obtained in variational form by Chen and M\"orters in~\cite{chen-morters}. 
Thus, it is given in terms of an operator $U_h$ by 
\[
\I_{CM}=\inf\acc{ ||h||_2:\ h\ge 0,\ \text{and}\ ||U_h||\ge 1},
\]
where
\be{rate-CM}
U_h(f)=\sqrt{e^{h}-1}\pare{G_d-\delta_0}*(f \sqrt{e^{h}-1}),
\ee
and $\delta_0$ is the delta function at 0, and for $z\in \Z^d$,
$f*g(z)=\sum_{y\in \Z^d} f(y)g(z-y)$.
\ep
\br{rem-limit}
We remark that Chen and M\"orters' proof produces (and relies on)
a finite volume version of their large deviation principle. 
Namely, for any finite subset $\Lambda\subset
\Z^d$,
\be{CM-finite}
\lim_{t\to\infty} \frac{1}{\sqrt t} \log \P
\pare{ \sum_{z\in \Lambda} l_\infty(z)\tilde l_\infty(z) \ge t}
=-2\I_{CM}(\Lambda),\quad\text{and}\quad \lim_{\Lambda \nearrow \Z^d}
\I_{CM}(\Lambda)=\I_{CM}.
\ee
However, \reff{CM-finite} gives no information on understanding which level sets of
the local times is responsible for the large deviation.
\er

Proposition~\ref{prop-limit} relies on the fact that
both the excess self-intersection local times and large intersection local times
are essentially realized on a finite region. The physical phenomenon behind
this observation is that the sites which contribute in
making $\{\langle l_\infty,\tilde l_\infty\rangle >t\}$ are those for which 
both $l_\infty(z)$ and $\tilde l_\infty(z)$ are
larger than $\sqrt{t}/A$ for some constant $A$. In other words,
define for $\xi>0$
\[
\D(\xi):=\{z\in \Z^d:\ l_\infty(z)\ge \xi\},\quad\text{and}\quad
\tilde \D(\xi):=\{z\in \Z^d:\ \tilde l_\infty(z)\ge \xi\}.
\]
Then, our next result reads as follows.

\bp{prop-CM} Assume the lattice has five or more dimension. Then
\be{ineq-CM}
\limsup_{A\to\infty} \limsup_{t\to\infty} \frac{1}{\sqrt t}\log
\P\pare{\sum_{z\not\in \D(\frac{\sqrt t}{A})\cap \tilde
\D(\frac{\sqrt t}{A})} l_\infty(z)\tilde l_\infty(z)>t}=-\infty.
\ee
\ep

Proposition~\ref{prop-CM} is based on the idea that 
$\langle l_\infty,\tilde l_\infty\rangle$ is not {\it critical} in the sense
that even when {\it weighting less} intersection local times, the strategy
remains the same. In other words, define for $q\le 2$
\be{append.1}
\zeta(q)=\sum_{z\in \Z^d} l_\infty(z)\tilde l_\infty^{q-1} (z).
\ee
Then, we have the following Lemma, interesting on its own.
\bl{lem-amine}
Assume that $d\ge 5$. For any $2\ge q>\frac{d}{d-2}$, there is $\kappa_q>0$ such that
\be{append.2}
\P\pare{\zeta(q)>t}\le \exp(-\kappa_q t^{\frac{1}{q}}).
\ee
\el

Our second task is to have some information on the geometry of the
intersection of two level sets of the local times of two
independent walks. For instance, we would like to address the following 
question: knowing the volume of the intersection of say $\acc{z:\ 
l_\infty(z)\sim n}$ and $\acc{z:\ \tilde l_\infty(z)\sim n}$ what can we say about
the capacity of the intersection set?

To state our results, we need some more notations.
For a finite subset $\Lambda\subset \Z^d$, we denote by $\S(\Lambda)$ the set
of permutations of the sites of $\Lambda$. Thus, $\gamma\in \S(\Lambda)$ is
written as $\gamma=(\gamma(1),\dots,\gamma(|\Lambda|))$, and we set $\gamma(0)=0$.
We define 
\[
T_\Lambda=\inf\acc{n>0:\ S_n\in \Lambda},\quad\text{and for }\quad z\in \Z^d,\quad 
H(z)=\inf\acc{n\ge 0:\ S_n=z}.
\]

Our key proposition is the following.
\bp{prop-level} Assume $d\ge 3$. There is a positive constant $c_d$ such that
for any $\Lambda$ finite subset of $\Z^d$, and
$\acc{n(z),\ z\in \Lambda}$ positive integers, we have
\be{intro.1}
P_0\pare{l_\infty(z)=n(z),\ \forall z\in \Lambda}\le (c_d \bar n)^{|\Lambda|}
(|\Lambda|!)^{d}e^{-\under{n} \capa(\Lambda)}\!\!\! \sum_{\gamma\in \S(\Lambda\bs\{0\})}
\prod_{i=1}^{|\Lambda\bs\{0\}|} P_{\gamma(i-1)}\pare{H(\gamma(i))<\infty},
\ee
where we set
\[
\bar n=\max_{z\in \Lambda} n(z),\quad \under{n}=\min_{z\in \Lambda} n(z),\quad
\text{and}\quad \capa(\Lambda)=\sum_{z\in \Lambda} P_z(T_\Lambda =\infty).
\]
\ep
\br{rem-level} 
Proposition~\ref{prop-level} is useful when dealing with {\it large level sets}
of the local times. In other words, 
we need $\under{n}\times \capa(\Lambda)\gg |\Lambda|$.
Since, we will see that $\capa(\Lambda)\ge c_d |\Lambda|^{1-2/d}$, the
range of applicability of \reff{intro.1} is
$\{(\under{n},|\Lambda|): \under{n}\gg |\Lambda|^{2/d}\}$.
We note that an alternative proof of Proposition~\ref{prop-CM} can be obtained using
Proposition~\ref{prop-level}, in dimension 6 or more.
\er
Let us now define more notations. For an integer $n$, let
\be{not.1}
\L(n)=\acc{z\in \Z^d: l_\infty(z)=n},\quad\text{and}\quad
\tilde \L(n)=\acc{z\in \Z^d: \tilde l_\infty(z)=n}.
\ee
With the same hypotheses as in Proposition~\ref{prop-level}, we have the following
useful Corollary.
\bc{cor-level}
Assume $d\ge 5$. There are positive constants $\kappa_d, C_d$ such that
for any positive integers $n,m,L$, we have
\be{intro.2}
P\pare{|\L(n)\cap\tilde \L(m)|\ge L}
\le  (C_d\ n\times m)^L (L!)^{2d}\quad \exp\pare{-\kappa_d\ (n+m)\ L^{1-\frac{2}{d}}}.
\ee
\ec
\br{rem-cor}
The inequality \reff{intro.2} is useful when $n+m\gg L^{2/d}$. A result in
the direction of $m$ and $n$ small, 
is Lemma 2 of \cite{KMSS} dealing with the intersection of
the ranges of two independent random walk, that we denote for any $n\in \N\cup\{
\infty\}$, $\RR_n=\acc{z:l_n(z)>0}$ and similarly for $\tilde \RR_n$. 
Assume that $d\ge 5$, and for any $\epsilon>0$, and $L$ large
enough
\be{kmss}
\exp\pare{-L^{1-\frac{2}{d}+\epsilon}}\le 
\P\pare{|\RR_\infty\cap \tilde \RR_\infty|\ge L}\le 
\exp\pare{-L^{1-\frac{2}{d}-\epsilon}}.
\ee
Finally, we remark that a large deviation principle has been established
by Bolthausen, den Hollander and van den Berg in~\cite{BBH},
for $\{ |\RR_n\cap \tilde \RR_n|\ge \xi n\}$ for large $n$ and $\xi>0$, or rather
for the continuous counterpart, that is the intersection of the volume of
two independent Wiener sausages. However, a large deviation principle
for $|\RR_\infty\cap \tilde \RR_\infty|$ is still open.
\er
A lower bound for $P(|\L(n)\cap\tilde \L(m)|\ge L)$ obtains following
the strategy proposed in~\cite{KMSS} in the proof of the lower bound
for their Theorem 4. Thus, we force the random walk $S_n$ (resp. $\tilde S_n$)
to make $n\times L^{1-\frac{2}{d}+\epsilon}$-returns to 0 
(resp. $m\times L^{1-\frac{2}{d}+\epsilon}$-returns to 0), making probable
that all sites of a ball of radius $r$, with $|B(0,r)|\ge L\ge |B(0,r/2)|$,
are visited $n$-times (resp. $m$-times). Thus, the following lower bound holds.
\bc{cor-kmss}[Corollary of Proposition 6 of \cite{KMSS}]
For any $\epsilon>0$, and for $L$ large enough, \cite{KMSS} prove that
\be{intro-LB}
P\pare{|\L(n)\cap\tilde \L(m)|\ge L}\ge \exp(-(n+m)L^{1-\frac{2}{d}+\epsilon}).
\ee
\ec
A direct consequence of Corollaries~\ref{cor-level} and~\ref{cor-kmss} is some
information about the geometry of the intersection of two {\it high} level
sets, knowing that the volume of the intersection is large. Thus, we formulate
it in the following way.
\bc{cor-geom}
Let $n_L$ and $m_L$ such that $\max(n_L,m_L)\gg L^{d/2}$
as $L$ goes to infinity, then for any $\epsilon>0$
\be{asymp-geometry}
\lim_{L\to\infty} 
P\pare{\capa\pare{\L(n_L)\cap\tilde \L(m_L)}\ge L^{1-\frac{2}{d}+\epsilon}
\big|\ |\L(n_L)\cap\tilde \L(m_L)|\ge L}=0.
\ee
\ec
The paper is organized as follows. In Section~\ref{sec-KMSS}, we prove
Lemma~\ref{lem-amine}, and then Proposition~\ref{prop-CM} as its corollary.
We then show Proposition~\ref{prop-limit} in Section~\ref{sec-limit}.
Then, we focus on Proposition~\ref{prop-level} in Section~\ref{sec-amine}, which
is technically the longest part of the paper. 
Finally, the proof of Corollary~\ref{cor-level} follows in Section~\ref{sec-cor}.
\section{Proofs of Proposition~\ref{prop-CM} and of Lemma~\ref{lem-amine}}
\label{sec-KMSS}
\subsection{Proof of Lemma~\ref{lem-amine}}
We assume $d\ge 5$. Lemma~\ref{lem-amine} can be thought
of as an interpolation inequality between 
Lemma 1 and Lemma 2 of \cite{KMSS}, whose proofs follow a classical pattern
(in statistical physics) of estimating all moments of $\zeta(q)$. 
This control is possible
since all quantities are expressed in terms of iterates of the Green's function,
whose asymptotics are well known (see for instance Theorem 1.5.4 of \cite{LAWLER}).

From \cite{KMSS}, it is enough that for a positive constant $C_q$,
we establish the following control on the moments
\be{append.3}
\forall n\in \N,\qquad \EE[\zeta(q)^n]\le  C_q^n (n!)^q.
\ee
First, noting that $q-1\le 1$, we use Jensen's inequality in the last inequality
\ba{append.4}
\EE[\zeta(q)^n]&\le & \sum_{z_1,\dots,z_n\in \Z^d} E_0\cro{
\prod_{i=1}^n l_\infty(z_i)}  E_0\cro{\prod_{i=1}^n l_\infty(z_i)^{q-1}}\cr
&\le & 
\sum_{z_1,\dots,z_n\in \Z^d} \pare{E_0\cro{\prod_{i=1}^n l_\infty(z_i)}}^{q}
\ea
If $\S_n$ is the set of permutation of $\acc{1,\dots,n}$ (with
the convention that for $\pi\in \S_n$, $\pi(0)=0$)
we have,
\be{append.5}
\begin{split}
E\cro{\prod_{i=1}^nl_\infty(z_i)}= &
\sum_{s_1,\dots,s_n\in \N} P_0(S_{s_i}=z_i,\ \forall i=1,\dots,n)\\
\le& \sum_{\pi\in \S_n} \sum_{s_1\le s_2\le\dots\le s_n\in \N}
P_0(S_{s_i}=z_{\pi(i)},\ \forall i=1,\dots,n)\\
\le & \sum_{\pi\in \S_n} \prod_{i=1}^n G_d\pare{z_{\pi(i-1)},z_{\pi(i)}} .
\end{split}
\ee
Now, by H\"older's inequality
\ba{append.9}
\sum_{z_1,\dots,z_n}\pare{\sum_{\pi\in \S_n} 
\prod_{i=1}^n G_d\pare{z_{\pi(i-1)},z_{\pi(i)}}}^q&\le & 
\sum_{z_1,\dots,z_n} (n!)^{q-1} \sum_{\pi\in \S_n} 
\prod_{i=1}^n G_d\pare{z_{\pi(i-1)},z_{\pi(i)}}^q\cr
&=& (n!)^q \sum_{z_1,\dots,z_n} \prod_{i=1}^n G_d\pare{z_{i-1},z_{i}}^q
\ea
Thus, classical estimates for the Green's function, \reff{append.9} implies that
\be{append.10}
\begin{split}
\sum_{z_1,\dots,z_n\in \Z^d} \pare{E_0\cro{\prod_{i=1}^n l_\infty(z_i)}}^{q}\le&
(n!)^q C^n \sum_{z_1,\dots,z_n} \prod_{i=1}^n (1+||z_i-z_{i-1}||)^{q(2-d)}\\
\le & (n!)^q C^n \pare{ \sum_{z\in \Z^d} (1+||z||)^{d(2-d)}}^n.
\end{split}
\ee
Thus, when $d\ge 5$ and $q>\frac{d}{d-2}$, we have a constant $C_q>0$ such
that
\be{append.7}
\EE[\zeta(q)^n]\le C_q^n(n!)^q 
\ee
The Lemma follows now by routine consideration.
\subsection{Proof of the Proposition~\ref{prop-CM}}
This follows easily from Lemma~\ref{lem-amine}. Indeed, for $q<2$
\be{append.20}
\begin{split}
\acc{\sum_{z\not\in\D(\frac{\sqrt t}{A})\cap \tilde \D(\frac{\sqrt t}{A})}
l_\infty(z)\tilde l_\infty(z)>t}\subset&
\acc{\sum_{l_\infty(z)\le \frac{\sqrt t}{A}
} l_\infty(z)^{q-1}\tilde l_\infty(z)> \frac{t}{2}\pare{\frac{A}{\sqrt t}}^{2-q}}\\
&\cup \acc{\sum_{\tilde l_\infty(z)\le \frac{\sqrt t}{A}
} l_\infty(z)\tilde l_\infty(z)^{q-1}> \frac{t}{2}\pare{\frac{A}{\sqrt t}}^{2-q}}.
\end{split}
\ee
Then, since $1>\frac{2-q}{2}$, Lemma~\ref{lem-amine} applied to \reff{append.10}
implies that for large $t$
\be{append.11}
\P\pare{\sum_{z\not\in \D(\frac{\sqrt t}{A})\cap \tilde \D(\frac{\sqrt t}{A})} 
l_\infty(z)\tilde l_\infty(z)>t}\le \exp\pare{ -\kappa_d A^{\frac{2-q}{q}}
t^{1/2}},\quad\text{since }\quad \frac{1}{q}(1-\frac{2-q}{2})=\frac{1}{2}.
\ee
\section{Proof of Proposition~\ref{prop-limit}}\label{sec-limit}
The proof of Proposition~\ref{prop-limit} relies on~\cite{A07}, and
we first recall some of its key steps we need here.
Our main result in \cite{A07} (called there Theorem 1.1) is that 
there is a positive constant $\I(2)$, such that for $\xi>0$ the following
limit exists
\be{old.1}
\lim_{n\to\infty} \frac{1}{\sqrt{n}}\log P_0\pare{
\bra{l_n,l_n}-E_0\cro{\bra{l_n,l_n}}> n\xi}= -\I(2) \sqrt{\xi}.
\ee
To prove \reff{old.1}, we observed that in order to produce
an excess self-intersection local times, the walk spends a time
of order $\sqrt n$ over a finite number of sites. Thus, for some
constant $A$, for any $\epsilon>0$, $n$ large enough and a constant
$\gamma$ depending on $A$
\be{old.2}
P_0\pare{ \bra{l_n,l_n}-E_0\cro{\bra{l_n,l_n}}> n(1+\epsilon)}\le 
n^\gamma P_0\pare{ \sum_{z\in \D_n({\sqrt n}/A)}l_\infty^2(z)> n},
\ee
where $\D_n(\xi)=\acc{z\in \Z^d:\ l_{n}(z)>\xi}$. Proposition 6.6 of \cite{A07}
allows us to relocate $\D_n({\sqrt n}/A)$ into a finite region $\Lambda_\epsilon$ of
$\Z^d$, whose diameter is independent of $n$, so that for $c>0$
\be{old.3}
P\pare{ \sum_{\D_n(\frac{\sqrt n}{A})}l_\infty^2(z)> n} 
\le e^{c\epsilon{\sqrt n}} 
P\pare{ \sum_{z\in\Lambda_\epsilon} l _\infty^2(z)> n,
\Lambda_\epsilon\subset\D_\infty(\frac{{\sqrt n}}{A}),\ l_\infty(0)
=\max_{\Lambda_\epsilon} l_\infty}.
\ee
Then, visiting ${\sqrt n}$-times each site of $\Lambda_\epsilon$ requires only
a total time of order ${\sqrt n}$, and Proposition 7.1 of \cite{A07} yields 
an integer sequence $\acc{k_n(z),z\in \Lambda_\epsilon}$ with
\[
\sum_{z\in \Lambda_\epsilon} k_n(z)^2\ge n,\quad k_n(z)\le A{\sqrt n},
\]
and, constants $\gamma,\alpha_0$ such that for $\alpha>\alpha_0$
\be{old.4}
P\pare{ \sum_{z\in\Lambda_\epsilon} l_\infty^2(z)> n,
\ l_\infty(0)=\max_{\Lambda_\epsilon} l_\infty}\le
n^\gamma P\pare{ l_{\lfloor \alpha\sqrt n\rfloor}|_{\Lambda_\epsilon}
=k_n|_{\Lambda_\epsilon}, \B(\lfloor \alpha\sqrt n\rfloor)},
\ee
where $\B(m)=\acc{S_m=0,l_m(0)=\max l_m(z)}$. On the other hand, note also
the obvious lower bound: for $\Lambda_\epsilon\subset \Lambda$
\be{old.5}
P\pare{ \sum_{\Lambda} l_\infty^2(z)> n}
\ge P\pare{ l_{\lfloor \alpha\sqrt n\rfloor}|_{\Lambda_\epsilon}=k_n|_{\Lambda_\epsilon},
\B(\lfloor \alpha\sqrt n\rfloor)}.
\ee
The subadditive argument treats the right hand side of \reff{old.5},
and yields also
\be{old.6}
\lim_{\Lambda\to\Z^d}\lim_{n\to\infty} \frac{1}{\sqrt n} \log P\pare{
\sum_{z\in\Lambda} l_\infty^2(z)> n} =-\I(2).
\ee
Now, we proceed with the link with intersection local times.
First, as mentioned in Remark~\ref{rem-limit},
Chen and M\"orters prove also that for any finite $\Lambda\subset \Z^d$
\[
\lim_{n\to\infty}\frac{1}{n^{1/2}}\log P\pare{ \bra{\ind_{\Lambda} l_\infty,
\tilde l_\infty}>n}=-2I_{CM}(\Lambda),
\]
with $I_{CM}(\Lambda)$ converging to $I_{CM}$ as $\Lambda$ increases to cover $\Z^d$.
The important feature is that for any fixed $\epsilon>0$,
we can fix a finite $\Lambda$ subset of $\Z^d$ such that $|I_{CM}(\Lambda)-I_{CM}|\le 
\epsilon$.
Note now that by Cauchy-Schwarz' inequality, and for finite set $\Lambda$
\be{old.7}
\sum_{z\in \Lambda} l_{\infty}(z)\tilde l_{\infty}(z)\le
\pare{\sum_{\Lambda} l_{\infty}^2(z)}^{\frac{1}{2}}
\pare{\sum_{\Lambda} \tilde l_{\infty}^2(z)}^{\frac{1}{2}}.
\ee
Inequalities \reff{old.6} and \reff{old.7} imply by routine consideration that
\be{old.8}
\limsup_{\Lambda\to\\Z^d} \limsup_{n\to\infty}  
\frac{1}{\sqrt n} \log P\pare{\bra{1_\Lambda l_{\infty},\tilde l_\infty}>n}
\le -\I(2)\inf_{\alpha>0}\acc{{\sqrt \alpha}+\frac{1}{{\sqrt \alpha}}}=-2\I(2).
\ee
Note also the obvious lower bound for $\Lambda_\epsilon\subset \Lambda$
\be{old.9}
P\pare{ \sum_{z\in \Lambda}l_\infty(z)\tilde l_\infty(z) > n}
\ge  P\pare{ l_{\lfloor \alpha\sqrt n\rfloor}|_{\Lambda_\epsilon}=
k_n|_{\Lambda_\epsilon}, \B(\lfloor \alpha\sqrt n\rfloor)}^2.
\ee
Since $\epsilon$ is arbitrary, we obtain
\be{old.10}
\liminf_{n\to\infty}  \frac{1}{\sqrt n} \log P\pare{\bra{l_{\infty},\tilde l_\infty}>n}
\ge -2\I(2).
\ee
\reff{old.8} and \reff{old.10} conclude the proof of Proposition~\ref{prop-limit}.
\section{Proof of Proposition~\ref{prop-level}}\label{sec-amine}
Site $0$ of $\Z^d$ plays a special r\^ole, since the walk start from 0.
We assume that $0\not\in \Lambda$, and omit to write the trivial changes
brought by the case $0\in \Lambda$. 

Our first step is to decompose paths in $\acc{l_\infty(z)=n(z),\ \forall z\in \Lambda}$ 
in terms of their sequence of crossed oriented edges.
Let $t=\sum_{z\in \Lambda} n(z)$, the total time the walk spends on $\Lambda$,
and let our state space be
\[
\Omega_t=\acc{\sigma\in (\{0\}\times \Lambda)\times(\Lambda\times\Lambda)^{t-1}
\times(\Lambda\times \{\infty\}):
\sigma_2(i)=\sigma_1(i+1),\forall i=1,\dots,t}.
\]
For an (oriented) edge $e\in \Lambda\times\Lambda$ (resp. for a site $z\in \Lambda$),
we denote by $l_\infty(e;.)$ (resp. $l_\infty(z;.)$) the variable on $\Omega_t$ counting
the crossings of $e$ (resp. the visits of $z$). Thus
\be{proof.1}
l_\infty(e;\sigma)=\sum_{s=1}^{t}\ind\acc{\sigma(s)=e},\quad\text{and}\quad
l_\infty(z;\sigma)=\sum_{s=1}^{t}\ind\acc{\sigma_2(s)=z}.
\ee
Now, to write concisely the decomposition we alluded to, we introduce
the following handy notation: for an edge $e=(x,y)$,
we set $q(e)=q(x,y)=P_x(T_\Lambda<\infty,\ S_{T_\Lambda}=y)$, 
\be{proof.4}
P_0\pare{l_\infty(z)=n(z),\ \forall z\in \Lambda}
=\!\!\sum_{\sigma\in \Omega_t}\!\!
\ind\acc{l_\infty(z;\sigma)=n(z),\ \forall z\in \Lambda}
\pare{ \prod_{i=1}^{t} q(\sigma(i))} P_{\sigma_2(t)}(T_\Lambda=\infty).
\ee

Now, we interpret $\sigma\in \Omega_t$ by specifying the fate of each site
of $\Lambda$. We associate, to each $z\in \Lambda$, a vector with $n(z)$ entries:
the $k$-th entry is the endpoint of the $k$-th visited edge with starting vertex $z$
(and if $z$ is the last visited site of $\Lambda$, then the $n(z)$-entry is $\infty$).
To identify fully the circuit, we need to specify $S_{T_\Lambda}$,
the first visited site of $\Lambda$ which we associate with 0. The family of vectors
(with the 1-vector associate with 0) is equivalent to $\sigma$, and specify a circuit
over $\Lambda$. Note that for any integer sequence 
$\acc{\E(e),\ e\in \Lambda^2}$, we have the multinomial
domination
\be{proof.3}
|\acc{\sigma\in \Omega_t:\ l_\infty(e;\sigma)=\E(e),\ \forall \ e\in \Lambda^2}|\le
\prod_{z\in \Lambda} \frac{\pare{\sum_{x\in \Lambda} \E(z,x)}!}
{\prod_{x\in \Lambda} \E(z,x)!}.
\ee

Since the event of interest is given in terms of visits of sites
of $\Lambda$, we need to express $l_\infty(z;\sigma)$ in terms 
of $\acc{l_\infty(e;\sigma),\ e\in \Lambda^2}$. 
The number of visits of $z\in \Lambda$ is essentially the
total crossings of all edges with initial vertex $z$. That is, for $\sigma\in \Omega_t$
\be{proof.2}
l_\infty(z;\sigma)=\sum_{x\in \Lambda} l_\infty((z,x);\sigma)+\ind\acc{z=\sigma_2(t)}.
\ee
For $z_1,z_t\in \Lambda$, let $\H(z_1,z_t)$ be the 
images of $\acc{l_\infty(e,\sigma),\ e\in \Lambda^2}$ as $\sigma$ spans $\Omega_t$ with
$\sigma_2(1)=z_1,\sigma_2(t)=z_t$, under the additional constraint that
for $\E\in \H(z_1,z_t)$, and
all $z\in \Lambda$, $\sum_{x\in \Lambda} \E(z,x)=n(z)-1\{z=z_t\}$.

We can now express \reff{proof.4} into a sum over edge crossings.
We need however to distinguish the first and last 
points in $\Lambda$, say $z_1$ and $z_t$ respectively. 
Thus,
\ba{proof.5}
P_0\pare{l_\infty(z)=n(z),\ \forall z\in \Lambda}
&\!\!\!\le\!\!\! &\sum_{z_1,z_t\in \Lambda} q(0,z_1)\sum_{\E\in \H(z_1,z_t)}
\prod_{e\in \Lambda^2} q(e)^{\E(e)}
\sum_{\sigma\in \Omega_t} 
\ind\acc{l_\infty(\bullet;\sigma)=\E(\bullet)}\cr
&\!\!\!\le\!\!\! &
\sum_{z_1,z_t\in \Lambda} q(0,z_1)\sum_{\E\in \H(z_1,z_t)}
\prod_{e\in \Lambda^2} q(e)^{\E(e)}
\prod_{z\in \Lambda} 
\frac{\pare{\sum_{x\in \Lambda} \E(z,x)}!}{\prod_{x\in \Lambda} \E(z,x)!}.
\ea
The proof proceeds now in two steps which we have gathered in two distinct
sections for the ease of reading.
\begin{itemize}
\item[Step 1:]Given $z_1,z_t\in \Lambda$ and 
the edge-occupation trace of path $\E=\acc{\E(e),e\in \Lambda^2}\in \H(z_1,z_t)$, 
we wish to extract
a path visiting each vertex of $\Lambda$ at least once and whose total (graph) length
grows like $|\Lambda|$ (and not $|\Lambda|^2$). We overcome this difficulty
by building a self-avoiding path covering $\Lambda$ whose set
of crossed-edges is not a subset of $\E$, but for which we have some domination
of the product over edge-distances (see \reff{tree.1}).
This is done in Section~\ref{sec-trail}.
\item[Step 2:]Extract the contribution of the self-avoiding path in \reff{proof.5}, and
re-adjust the sum in \reff{proof.5} so that a true multinomial distribution appears, as
well as the capacity of $\Lambda$. This is done in Section~\ref{sec-capacity}.
\end{itemize}
\subsection{Extracting a Self-Avoiding path}\label{sec-trail}
We fix $z_1,z_t\in \Lambda$, and $\E\in \H(z_1,z_t)$. Let $\gamma$ be a path
giving rise to the edge-occupation numbers $\E$. Our aim in this section is
to extract a self-avoiding path over $\Lambda$ which {\it uses} the edges of $\E$.
The self-avoiding path we eventually build is not a sub-path of $\gamma$, though
it starts with $z_1$. To describe in details the construction, we first set notations.
A self-avoiding path over $\Lambda$, which we call for simplicity
a {\it $\Lambda$-trail}, is an element of $(\Lambda^2)^{|\Lambda|-1}$,
say $\acc{e_1,\dots,e_{|\Lambda|-1}}$ such that for $i=1,\dots,|\Lambda|-2$, 
the end-vertex
of $e_i$ is the starting-vertex of $e_{i+1}$, and so that no vertex is used
twice.  The set of loops of $\Lambda$ is denoted 
$\Delta_\Lambda=\acc{(z,z): z\in \Lambda}$.
Now, for each $z_1,z_t\in \Lambda$, and $\E\in \H(z_1,z_t)$, we call a $\E$-stock, an
element $\S\in \N^{\Lambda\times\Lambda}$ with $\S\le \E$ for the natural order
in $\N^{\Lambda\times\Lambda}$, with $\S|_\Delta\equiv 0$, $\sum_{e\in 
\Lambda\times\Lambda} \S(e)=|\Lambda|-1$, and such that for each
vertex $z\in\Lambda$, there is at most one oriented edge $e\in \acc{z}\times
\Lambda$ such that $\S(e)>0$. Let us denote by $(\S)$ this latter property. 
Note that contrary to the
edge-occupation number $\E$, a $\E$-stock may not correspond to a path.

We show in this section that for each $z_1,z_t\in \Lambda$, and $\E\in \H(z_1,z_t)$,
there is a $\Lambda$-trail $\T$, with initial vertex $z_1$,
and a $\E$-stock $\S$ such that if for edge $e=(z,z')$ we denote by
$d(e)=||z-z'||$, then
\be{tree.1}
|\Lambda|!\prod_{e\in \Lambda^2} d(e)^{\S(e)}\ge \prod_{e\in \T} d(e).
\ee
\subsubsection{Coalescing $\Lambda$}
Let $N=|\Lambda|$, and call $\Lambda_N=\Lambda$ to emphasize that $\Lambda$ contains
$N$ vertices. We coalesce $\Lambda$ through a sequence of graphs $\Lambda_{N-1},
\dots,\Lambda_1$ such that for $k=1,\dots,N$ the graph $\Lambda_k$ has
exactly $k$ vertices. Then, we show \reff{tree.1} by induction on $\Lambda_k$.
We first describe how to obtain $\Lambda_{N-1}$. Choose a pair of vertices,
say $e_N=(a_N,\tilde a_N)$ of $\Lambda$ which minimize the euclidean distance.
$\Lambda_{N-1}$ is obtained from $\Lambda_N$ by suppressing $a_N,\tilde a_N$
but adding one new vertex, say $\psi_{N-1}$, which we also think of as
a cluster containing $a_N$ and $\tilde a_N$. We also define a pseudo-metric
$d_{N-1}:\Lambda_{N-1}^2\to\R^+$ by
\[
\forall z,z'\in \Lambda_{N-1}\bs \{\psi_{N-1}\},\qquad
d_{N-1}(z,z')=d(z,z'),
\]
and 
\[
d_{N-1}(z,\psi_{N-1})=d_{N-1}(\psi_{N-1},z)=\min\pare{ d(z,a_N),d(z,\tilde a_N)}.
\]
Note that $d_{N-1}$ fails to satisfy the triangle inequality. Now, we can
associate with $\gamma$ and 
$\Lambda_{N-1}$ an edge-occupation number $\E_{N-1}$. In other words,
from path $\gamma$, we count edge-crossings as if $a_N$ and $\tilde a_N$ were
indistinguishable, and obtain a path on $\Lambda_{N-1}$ which we call
$\gamma_{N-1}$. Thus, for $z\in \Lambda_{N-1}\bs \{\psi_{N-1}\}$
\be{edge-rec}
\E_{N-1}(z,\psi_{N-1})=\E(z,a_N)+\E(z,\tilde a_N),\quad 
\E_{N-1}(\psi_{N-1},z)=\E(a_N,z)+\E(\tilde a_N,z),
\ee
and, $\E_{N-1}(\psi_{N-1},\psi_{N-1})=\E(a_N,\tilde a_N)+\E(\tilde a_N,a_N)+
\E(a_N,a_N)+\E(\tilde a_N,\tilde a_N)$. We proceed by induction until we 
reach $\Lambda_1$ with one vertex, which we think of as a cluster of $N$ vertices
of $\Lambda$.

Note that for $k=2,\dots,N$, and $z,z'\in \Lambda_k$, we keep in mind
that $z,z'$ can be thought of as clusters in $\Lambda_N$, and that
\be{pseudo-d}
d_k(z,z')=\min\acc{d(x,x'):\ x,x'\in \Lambda_N,\ x\in z,\ x'\in z'}.
\ee
Also, we call $\gamma_k$ the path obtained from $\gamma$ on $\Lambda_k$
by identifying all vertices (of $\Lambda$) belonging to a single cluster (i.e. a vertex
of $\Lambda_k$).
\subsubsection{Proof of \reff{tree.1} by induction}
We want that $z_1$ be the starting point of the trail. Note that $\Lambda_2$
has two vertices $\acc{a_2,\tilde a_2}$, and since all sites of $\Lambda$
are visited, we have necessarily $\E_2(a_2,\tilde a_2)>0$ when $z_1\in a_2$
and $\E_2(\tilde a_2,a_2)>0$ when $z_1\in \tilde a_2$, so we can
define an $\E_2$-block $\S_2(a_2,\tilde a_2)=1$ (resp. $\S_2(\tilde a_2,a_2)=1$),
and the $\Lambda_2$-trail
$\T_2=\acc{(a_2,\tilde a_2)}$ (resp. $\T_2=\acc{(\tilde a_2,a_2)}$),
when $z_1$ belongs to the first (resp. second) cluster.

Assume now that we have an $\E_{k-1}$-stock $\S_{k-1}$, and a $\Lambda_{k-1}$-trail
$\T_{k-1}$ so that
\be{ind-k}
\frac{N!}{(N-(k-2))!}\prod_{e\in \Lambda_{k-1}^2} d_{k-1}(e)^{\S_{k-1}(e)}
\ge \prod_{e\in \T_{k-1}} d_{k-1}(e).
\ee
Recall that $e_k=(a_k,\tilde a_k)$ has been coalesced to produce vertex $\psi_{k-1}$
of $\Lambda_{k-1}$. Note that trail $\T_{k-1}$ crosses $\psi_{k-1}$ only
once. Also, it is part of our induction hypothesis to assume that the first
crossed vertex of $\T_{k-1}$ when seen as a cluster, contains $z_1$.

Thus, assume first that $\psi_{k-1}$ is not the first vertex of the trail $\T_{k-1}$.
Let $b,b'\in \Lambda_{k-1}$ such that
$(b,\psi_{k-1})$ and $(\psi_{k-1},b')\in \T_{k-1}$. 

By construction, $\E_k$ satisfies for any $z\in \Lambda_k$
\be{tree.3}
\E_k(z,a_k)+\E_k(z,\tilde a_k)=\E_{k-1}(z,\psi_{k-1})\ge \S_{k-1}(z,\psi_{k-1}),
\ee
and,
\[
\E_k(a_k,z)+\E_k(\tilde a_k,z)=\E_{k-1}(\psi_{k-1},z)\ge \S_{k-1}(\psi_{k-1},z).
\]
Thus, we can define $\tilde \S_{k}$ on $\Lambda_k^2$ such that for $z\in \Lambda_k
\cap \Lambda_{k-1}$
\[
\tilde \S_k(z,a_k)+\tilde \S_k(z,\tilde a_k)=\S_{k-1}(z,\psi_{k-1}),\quad
\tilde \S_k(z,a_k)\le \E_k(z,a_k),\quad\text{and}\ 
\tilde \S_k(z,\tilde a_k)\le \E_k(z,\tilde a_k).
\]
and similarly for $\tilde \S_k(a_k,z)$ and $\tilde \S_k(\tilde a_k,z)$.
Also, for $z,z'\in \Lambda_k\bs\acc{a_k,\tilde a_k}$, $
\tilde S_k(z,z')=\S_{k-1}(z,z')$.

Note that in path $\gamma_k$, there is at least one crossing of each
vertex of $\Lambda_k$. Thus,
\be{tree.10}
\sum_{e\in \Lambda_k\bs\Delta_k}\E_k(e)\ge k-1, \quad\text{where}\quad
\Delta_k:=\Delta_{\Lambda_k}.
\ee
Now, by the induction hypothesis, $\S_{k-1}$ satisfies $\sum_e \S_{k-1}(e)=k-2$ and 
$(\S)$. The simple observation which guarantees property $(\S)$ to
propagate through induction is that there cannot be two vertices of $\Lambda_k$,
say $x$ and $y$, such that 
\[
\forall z\in \Lambda_k\bs\acc{x} \quad\E_k(x,z)=0,\quad\text{and}\quad 
\forall z\in \Lambda_k\bs\acc{y}\quad \E_k(y,z)=0.
\]
The reason is that when looking at $\gamma_k$, at most only one of $x$ or $y$
can be the last visit in $\Lambda_k$ (and in this case, the last vertex points to no
other vertex of $\Lambda_k$). We consider now two cases. Assume first, that 
\[
\forall z\in \Lambda_k\bs\acc{a_k}\quad  \tilde S_k(a_k,z)=0,\quad\text{and}\quad
\forall z\in \Lambda_k\bs\acc{\tilde a_k}\quad \tilde S_k(\tilde a_k,z)=0.
\]
Then, we just saw that there is $z^*\in \Lambda_k$ such that
$\E_k(a_k,z^*)>0$ and $z^*\not= a_k$, or $\E_k(\tilde a_k,z^*)>0$ and 
$z^*\not= \tilde a_k$. In the case
$\E_k(a_k,z^*)>0$ (resp. $\E_k(\tilde a_k,z^*)>0$),
 we choose $e_k^*=(a_k,z^*)$ (resp. $e_k^*=(\tilde a_k,z^*)$), 
and $\S_k=\tilde\S_{k}+1_{e_k^*}$. Secondly, assume that for $z^*\in \Lambda_k$,
we have $\tilde S_k(a_k,z^*)>0$ and necessarily $\tilde S_k(\tilde a_k,z)=0$
for all $z\in \Lambda_k$ (the other case 
$\tilde S_k(\tilde a_k,z^*)>0$ is similar). By the induction hypothesis,
there is another vertex of $\Lambda_k$, say $a^*$ such that for any
$z\in \Lambda_k$, we have $\tilde S_k(a^*,z)=0$.
We then proceed as in the first case to deduce that there is an edge $e_k^*$
incident to either $\tilde a_k$ or $a^*$ such that $\E_k(e_k^*)>0$, and $e_k^*$
is not a loop.

Thus, there
is an edge $e_k^*\in \Lambda_k^2\bs\Delta_k$ such that $\E_k(e_k^*)>\tilde
\S_k(e_k^*)$, where $\tilde \S_k$ is built from $\S_{k-1}$ as in the previous case.
We thus set $S_k=\tilde S_k+1_{e_k^*}$, and note that by definition
\be{tree.11}
d_k(a_k,\tilde a_k)\le d_k(e_k^*).
\ee
For the ease of reading we distinguish two cases.

\noindent{ \underline{Case 1}: 
$d_{k-1}(b,\psi_{k-1})=d_{k}(b,a_k)$ and $d_{k-1}(\psi_{k-1},b')=d_k(\tilde a_k,b')$} 

The trail $\T_k$ is the same as trail $\T_{k-1}$
but with edges $\acc{(b,\psi_{k-1}),(\psi_{k-1},b')}$ replaced by edges
$\acc{ (b,a_k),(a_k,\tilde a_k),(\tilde a_k,b')}$. Thus, we have
\be{tree.4}
d_k(e_k^*)\times d_{k-1}(b,\psi_{k-1})d_{k-1}(\psi_{k-1},b')\ge 
d_k(b,a_k)d_k (a_k,\tilde a_k)d_k (\tilde a_k,b').
\ee
The same reasoning (with obvious changes of symbols) would hold if
$d_{k-1}(b,\psi_{k-1})=d_k(\tilde a_k,b)$ and $d_{k-1}(\psi_{k-1},b')=
d_k(a_k,b')$.

\noindent{ \underline{Case 2}: 
$d_{k-1}(b,\psi_{k-1})=d_{k}(b,a_k)$ ,and $d_{k-1}(\psi_{k-1},b')= d_k(a_k,b')$.}

We have to think of $a_k,\tilde a_k, b'$ as clusters, and let
$x,x'\in a_k$, $y,y'\in \tilde a_k$ and $z,z'\in b'$ be vertices of $\Lambda_N$
such that
\be{tree.5}
d(x',z)=d_k(a_k,b'),\quad d(x,y')=d_k(a_k,\tilde a_k),\quad\text{and}\quad
d(y,z')=d_k(\tilde a_k,b').
\ee
Since $x,x'$ belong to the same cluster, our construction implies that there 
is a sequence $\acc{x_0,x_1,\dots,x_{k_1}}$ vertices of $a_k$ with
\be{tree.6}
x_0=x,\quad x_{k_1}=x',\quad\text{and}\quad \forall i=1,\dots,k_1\quad
d(x_{i-1},x_i)\le d_k(a_k,\tilde a_k).
\ee
Similarly, we consider $\acc{y_0,y_1,\dots,y_{k_2}}$ in $\tilde a_k$ joining
$y$ and $y'$ and $\acc{z_0,z_1,\dots,z_{k_3}}$ in $b'$ joining $z$ and $z'$.
The important observation is that $k_1+k_2+k_3\le N-k$. Thus, by the triangle 
inequality
\be{tree.7}
\begin{split}
d_k(\tilde a_k,b')=d(y,z') &\le \sum_{i=1}^{k_2}d(y_{i-1},y_i)+d(y',x)+
\sum_{i=1}^{k_1}d(x_{i-1},x_i)+d(x',z)+\sum_{i=1}^{k_3}d(z_{i-1},z_i)\\
&\le (N-k+1)d_k(a_k,\tilde a_k)+d_k(a_k,b')\\
&\le (N-k+2) d_k(a_k,b')
\end{split}
\ee
Using \reff{tree.7}, we have
\be{tree.8}
(N-k+2)d_k(e_k^*) d_{k-1}(b,\psi_{k-1})d_{k-1}(\psi_{k-1},b')\ge
d_{k}(b,a_k) d_k (a_k,\tilde a_k)d_k (\tilde a_k,b').
\ee
Thus, with the same trail as in the previous case, we have
from \reff{ind-k} and \reff{tree.8}
\be{tree.9}
\frac{N!}{(N-k+1)!} \prod_{e\in \Lambda_k^2} d_k(e)^{S_k(e)}\ge
\prod_{e\in \T_k} d_k(e).
\ee

Finally, we explain how to manage so that $z_1$ remains in the first cluster
of $\T_k$. Assume that $\psi_{k-1}$ contains $z_1$ (and necessarily $\psi_{k-1}$
is the first vertex of $\T_{k-1}$). If $z_1\in a_k$, then add a fictitious
vertex $b^{\dagger}$ in $\Lambda_k$, and an edge $(b^{\dagger},a_k)$, and
proceed as before with $b^{\dagger}$ in the r\^ole of $b$. However, if
$z_1\in \tilde a_k$, the fictitious edge is $(b^{\dagger},\tilde a_k)$, and
$\tilde a_k$ replaces $a_k$ in the previous constructions.

\subsection{About Outer Capacity}\label{sec-capacity}
For $z_1,z_t\in \Lambda$  and $\E\in \H(z_1,z_t)$, let $\S$ be the $\E$-stock built
in Section~\ref{sec-trail}.
The first step is to subtract edges of $\S$ from $\E$,
and replace them with loops of $\Delta$.
Also, the number of edges of $\S$ incident with $z\in \Lambda$, is denoted
$\leaf(z,\S)$, that is $\leaf(z,\S)=\sum_{z'} \S(z,z')\in \{0,1\}$. 
Thus, we define for $e\not\in \Delta$
\be{capa.1}
\E_{\S}(e)=\E(e)-\S(e)\ge 0,
\ee
and,
\be{capa.2}
\E_{\S}(z,z)=\E(z,z)+\leaf(z,\S) \quad(\text{recall that } \S(z,z)=0).
\ee
Note that we might have $\E_{\S}\not\in \H(z_1,z_t)$, but
at least each vertex $z$ of $\Lambda$ occurs $n(z)$ times, in the sense that
\be{capa.3}
\forall z\in \Lambda, \qquad \sum_{z'\in \Lambda} \E_{\S}(z,z')=\sum_{z'\in \Lambda} \E(z,z').
\ee
The important (and simple) observation is that for $\E\in \H(z_1,z_t)$
\be{capa.4}
\prod_{z\in \Lambda}\frac{\pare{\sum_{x\in \Lambda} \E(z,x)}!}{\prod_{x\in \Lambda} \E(z,x)!}
\le \pare{ \max_{z\in \Lambda} n(z)}^{|\Lambda|}
\prod_{z\in \Lambda}
\frac{\pare{\sum_{x\in \Lambda} \E_{\S}(z,x)}!}{\prod_{x\in \Lambda} \E_{\S}(z,x)!}
\ee
Indeed, in view of \reff{capa.3}, \reff{capa.4} is equivalent to showing that
\be{capa-heart}
\prod_{e\in \Delta_\Lambda} \frac{\E_{\S}(e)!}{\E(e)!}\le \bar n^{|\Lambda|}
\prod_{e\not\in \Delta_\Lambda} \frac{\E(e)!}{\E_{\S}(e)!}.
\ee
Now since $\E(e)\ge \E_{\S}(e)$ for $e\not\in\Delta_\Lambda$, 
Now, \reff{capa-heart} would follow from
\[
\prod_{z\in \Lambda} \frac{(\E(z,z)+\leaf(z,\S))!}{\E(z,z)!}\le \bar n^{|\Lambda|}.
\]
Note that $\E(z,z)+\leaf(z,\S)\le \sum_x \E(z,x)\le n(z)$, so that
\[
\frac{(\E(z,z)+\leaf(z,\S))!}{\E(z,z)!}\le n(z)^{\leaf(z,\S)},
\]
and since $\sum_z \leaf(z,\S)=\sum_e \S(e)= |\Lambda|-1$, 
we have for $e\in \Delta_\Lambda$
\be{capa-obvious}
\prod_{z\in \Lambda} \frac{(\E(z,z)+\leaf(z,\S))!}{\E(z,z)!}\le 
\prod_{z\in \Lambda} n(z)^{\leaf(z,\S)}\le (\bar n)^{|\Lambda|}.
\ee
Also, since the walk has probability $\frac{1}{2d+1}$ to stay still, we have
\be{capa.5}
q(e)\ge \frac{1}{2d+1},
\ee
so that using that $\sum_e \S(e)=|\Lambda|-1$ 
\be{capa.6}
\prod_{e\in \Lambda^2} q(e)^{\E(e)-\S(e)} \le (2d+1)^{|\Lambda|}
\prod_{e\in \Lambda^2} q(e)^{\E_{\S}(e)}.
\ee
Now, we call $\TT(z_1)$ the set of $\Lambda$-trails with initial vertex $z_1$,
and combine \reff{capa.4}, \reff{capa.6} into \reff{proof.5} to obtain
\ba{capa.7}
P_0\pare{l_\infty(z)=n(z),\ \forall z\in \Lambda}
&\le &\pare{(2d+1)\bar n}^{|\Lambda|}
\sum_{z_1,z_t\in \Lambda} q(0,z_1)\cr
\sum_{\T\in \TT(z_1)}
\sum_{\E\leftrightarrow \T}\!\!\!&&\!\!\!\prod_{e\in \Lambda^2} q(e)^{\S(e)}
\prod_{e\in \Lambda^2} q(e)^{\E_{\S}(e)}
\prod_{z\in \Lambda} \frac{\pare{\sum_{x\in \Lambda}\E_{\S}(z,x)}!}
{\prod_{x\in \Lambda}\E_{\S}(z,x)!} .
\ea
By $\E\leftrightarrow \T$ we mean that $\E\in \H(z_1,z_t)$
gives rise to trail $\T$, and by $\S$ we denote the $\E$-stock associated with $\E$. 

We make three observations.
\begin{enumerate}
\item For any fixed $z\in \Lambda$, 
\be{capa.8}
\prod_{x\in \Lambda} P_z(S_{T_\Lambda}=x)^{\E(z,x)}=
P_z(T_\Lambda<\infty)^{\sum_x \E(z,x)}\times 
\prod_{x\in \Lambda} p(z,x)^{\E(z,x)},
\ee
where, for $x,z\in \Lambda$ we defined
\be{capa.10}
p(z,x)=P_z(S_{T_\Lambda}=x|T_\Lambda<\infty),
\quad\text{wich satisfy}\quad 
\sum_{x\in \Lambda} p(z,x)=1.
\ee
Also, using that $1-x\le \exp(-x)$ (and with $\under{n}=\min(n(z))$)
\ba{capa.9}
\prod_{z\in \Lambda} P_z(T_\Lambda<\infty)^{n(z)-1\{z=z_t\}}\le&&
\exp\pare{-\sum_{z\in \Lambda} P_z(T_\Lambda=\infty)(n(z)-1\{z=z_t\})}\cr
\le &&\exp\pare{1-\under{n}\sum_{z\in \Lambda} P_z(T_\Lambda=\infty)}
=e^{1-\under{n}\ \capa(\Lambda)},
\ea
where $\capa(\Lambda)$ denotes the capacity of the finite set $\Lambda$, given
by $\capa(\Lambda)=\sum_{z\in \Lambda} P_z(T_\Lambda=\infty)$ (see Section 2.2
of~\cite{LAWLER}).
\item For any edge $e$ joining $z_1,z_2\in \Lambda$ 
\be{capa.13}
q(e)=P_{z_1}(S_{T_\Lambda}=z_2) \le  P_{z_1}(H(z_2)<\infty),
\ee
and well known asymptotics (see e.g. Theorem 1.5.4 of \cite{LAWLER}), say
that for $d\ge 3$, there are positive constants $c_d,
{\underline c_d}$ such that for any $z_1$ and $z_2$ distinct vertices of $\Z^d$
\be{capa.14}
{\underline c_d} ||z_1-z_2||^{2-d}\le P_{z_1}(H(z_2)<\infty)\le  c_d ||z_1-z_2||^{2-d}.
\ee
Thus, recalling that for $e=(z_1,z_2)$ we set 
$d(e):=||z_1-z_2||$, inequality \reff{tree.1} of Section~\ref{sec-trail}
yields for $\E$-block $\S$ and $\Lambda$-trail $\T$
\be{capa.15}
\prod_{e\in \Lambda^2} q(e)^{\S(e)}\le c_d^{|\Lambda|-1}
\pare{ \prod_{e\in \Lambda^2} d(e)^{\S(e)}}^{2-d}\le
c_d^{|\Lambda|-1}(|\Lambda|!)^{d-2}\pare{ \prod_{e\in \T} d(e)}^{2-d}.
\ee
\item Because of property $(\S)$, the transformation $\E\to\E_\S$ can send
at most $(|\Lambda|-1)^{|\Lambda|}$ edge-crossing configurations
to the same image. Indeed, for each vertex
we have $|\Lambda|-1$ possible edges (to the $|\Lambda|-1$ distinct vertices) 
which could have been changed into a self-loop. Also, $(|\Lambda|-1)^{|\Lambda|}\le
e^{|\Lambda|} |\Lambda|!$.
\end{enumerate}
Thus, with $c_d'=2e^2c_d(2d+1)$, \reff{capa.7} becomes
\ba{capa.12}
P_0\pare{l_\infty(z)=n(z),\ \forall z\in \Lambda}
&\le &\pare{c_d' \bar n}^{|\Lambda|}(|\Lambda|!)^{d-2} e^{-\under{n} \capa(\Lambda)}
\sum_{z_1,z_t\in \Lambda} q(0,z_1)\sum_{\T\in \TT(z_1)}\prod_{e\in \T} d(e)^{2-d}\cr
&&\qquad\times \sum_{\E\leftrightarrow \T}
\prod_{e\in \Lambda^2} p(e)^{\E_{\S}(e)}
\prod_{z\in \Lambda} \frac{\pare{\sum_{x\in \Lambda}\E_{\S}(z,x)}!}
{\prod_{x\in \Lambda}\E_{\S}(z,x)!} .
\ea
Thus, from \reff{capa.12}, we obtain after summing over $\E_{\S}$,
and taking into account the degeneracy explained in point 3 above,
\be{capa.11}
\begin{split}
P_0(l_\infty(z)=n(z),&\ \forall z\in \Lambda)
\le \pare{e c_d'\bar n}^{|\Lambda|}|\Lambda|(|\Lambda|!)^{d-1} 
\ e^{-\under{n} \capa(\Lambda)} 
\sum_{\gamma\in \S(\Lambda)}
\prod_{i=1}^{|\Lambda|} d\pare{\gamma(i),\gamma(i-1)}^{2-d}\\
\times \sum_{ \E\in \N^{\Lambda\times\Lambda}}&
\ind_{\acc{n(z)=\sum_{x} \E(z,x)+ 1\acc{z=z_t},\ \forall z \in \Lambda}}
\prod_{z\in \Lambda} \pare{ \frac{\pare{\sum_{x\in \Lambda} \E(z,x)}!}
{\prod_{x\in \Lambda} \E(z,x)!} \prod_{x\in \Lambda} p(z,x)^{\E(z,x)}}\\
\le &\pare{e c_d'\bar n}^{|\Lambda|} (|\Lambda|!)^{d}\ 
e^{-\under{n}\ \capa(\Lambda)} \sum_{\gamma\in \S(\Lambda)}
\prod_{i=1}^{|\Lambda|} d\pare{\gamma(i),\gamma(i-1)}^{2-d}.
\end{split}
\ee
Proposition~\ref{prop-level} follows now from the \reff{capa.11} and the lower
bound in \reff{capa.14}.
\section{Proof of Corollary~\ref{cor-level}}\label{sec-cor}
Assume that $d\ge 5$. First, note that
\be{cor.0}
\P\pare{|\D(n)\cap \tilde \D(m)|\ge L}\le 
\sum_{\Lambda\subset \Z^d: |\Lambda|=L}
P_0\pare{ l_\infty(z)=n,\ \forall z\in \Lambda}\times
P_0\pare{ l_\infty(z)=m\ \forall z\in \Lambda}.
\ee
Thus, if we set $\v(x)=||x||^{2-d}$ for $x\in \Z^d$, 
Proposition~\ref{prop-level} yields by Cauchy-Schwarz, 
\be{cor.01}
\begin{split}
\P\pare{|\D(n)\cap \tilde \D(m)|\ge L}\le& \sum_{\Lambda\subset \Z^d: |\Lambda|=L}
(c_d n\times m)^{|\Lambda|}
e^{-(n+m) \capa(\Lambda)} (|\Lambda|!)^{2d}\\
&\qquad \sum_{\gamma\in \S(\Lambda\bs\{0\})}
\prod_{i=1}^{|\Lambda\bs\{0\}|} \v\big(\gamma(i)-\gamma(i-1)\big)^{2},
\end{split}
\ee
We first show that there is a constant $\kappa_d>0$, such that
for any finite $\Lambda\subset \Z^d$, we have
\be{cor.1}
\capa(\Lambda)\ge \kappa_d |\Lambda|^{1-\frac{2}{d}}.
\ee
We use the variational characterisation of capacity (see the Appendix of ~\cite{FU}),
which says that if $G_d$ is the Green kernel, $\mu$ is a non-negative measure 
on $\Lambda$, and $c$ a positive constant
\be{cor.2}
\forall z\in \Lambda,\quad \sum_{z'\in \Lambda} G_d(z,z')\mu(z')\le c
\Longrightarrow \capa(\Lambda)\ge \frac{\sum_{z\in \Lambda} \mu(z)}{c}.
\ee
We have shown in the proof of Lemma 1.2 of \cite{AC04}, that 
there is $\kappa_d>0$ such that for any $z\in \Lambda$
\be{cor.3}
\sum_{z'\in \Lambda} G_d(z,z')\mu(z')\le \frac{1}{\kappa_d},\quad\text{with}\quad
\mu(z)=\frac{\ind_{\Lambda}(z)}{|\Lambda|^{2/d}}.
\ee
The desired bound \reff{cor.1} follows readily from \reff{cor.3}.

Now, note that (with the convention $\gamma(0)=z(0)=0$)
\be{cor.6}
\sum_{\Lambda: |\Lambda|=L}\sum_{\gamma\in \S(\Lambda)}
\prod_{i=1}^{|\Lambda|} \v\pare{\gamma(i)-\gamma(i-1)}^2=
\sum_{z_1,\dots,z_L\text{ distinct}} \prod_{i=1}^{L} \v\pare{z_i-z_{i-1}}^2\le
\pare{ \sum_{z\in \Z^d} \v(z)^2}^L.
\ee
Now, dimension $5$ or more enters in making the last series convergent, and
the result follows at once.

\end{document}